\documentstyle[12pt]{article}                            
\newtheorem{theorem}{Theorem}[section]
\newtheorem{definition}[theorem]{Definition}
\newtheorem{proposition}[theorem]{Proposition}
\newtheorem{corollary}[theorem]{Corollary}
\newtheorem{lemma}[theorem]{Lemma}

\def \proof {\noindent {\bf Proof.}\ \ }
\def \endproof {{\mbox{}\nolinebreak\hfill\rule{2mm}{2mm}\par\medbreak} }

\newcommand {\text}[1] {\mbox{\ \ #1 \ \ }}

\def \# {\vspace{0.5cm} \noindent}

\def \e {\varepsilon}
\def \d {\delta}
\def \E {{\bf E}}
\def \P {{\bf P}}
\def \R {{\bf R}}
\def \xj {\xi_j}
\def \tpx {T_p(X)}
\def \tpxx {T_{p^*}(X^*)}
\def \lin {l_\infty^n}
\def \ra {\rightarrow}
\def \sign {{\rm sign}}
\def \aconv {{\rm a.conv}}
\def \span {{\rm span}}
\def \vol {{\rm Vol}}
\def \dm {\left( \begin{array}{c}  |D| \\ m  \end{array}  \right)}

\begin{document}
\title {Absolutely representing systems, uniform smoothness, and type}
\author {R.Vershynin}
\maketitle

\begin{abstract}
Absolutely representing system (ARS) in a Banach space $X$ is a set
$D \subset X$ such that every vector $x$ in $X$ admits a 
representation by an absolutely convergent series
$x = \sum_i a_i x_i$ with $(a_i)$ reals and $(x_i) \subset D$.
We investigate some general properties of ARS. 
In particular, ARS in uniformly smooth and 
in B-convex Banach spaces are characterized 
via $\varepsilon$-nets of the unit balls.
Every ARS in a B-convex Banach space is quick, i.e. in the representation
above one can achieve $\|a_i x_i\| < cq^i\|x\|$, $i=1,2,\ldots$
for some constants $c>0$ and $q \in (0,1)$.
\end{abstract}

\section{Introduction}
%----------------------------------------------------------------------------

The concept of absolutely representing system (ARS) goes back to 
Banach and Mazur (\cite{B}, p.~109--110).

\begin{definition}
  A set $D$ in a Banach space $X$ is called {\sl absolutely
  representing system (ARS)} if for every $x \in X$ there are
  scalars $(a_i)$ and elements $(x_i) \subset D$ such that
  $$
    x=\sum_{i=1}^\infty {a_i x_i}    \text{and} 
    \sum_{i=1}^\infty {\|a_i x_i\|} < \infty. 
  $$
\end{definition} 

\noindent It can be observed (Section \ref{Gen}) that 
if $D$ is an ARS, then there exist a constant $c$ such that each 
$x \in X$ admits a representation 
$x=\sum_{i=1}^\infty {a_i x_i}$ with $\sum{\|a_i x_i\|} \le c\|x\|$. 
Then we call $D$ a "$c$-ARS".

For needs of complex analysis, ARS were defined also in locally
convex topological spaces \cite{K81}. In the theory of analytical 
functions such ARS happen to be a convenient tool: see \cite{K96}, 
\cite{G}, \cite{A}. Many results 
of general kind on ARS are obtained by Yu. Korobe\u\i nik and his
collaborators: see, for example, \cite{K81}, \cite{K86}, \cite{KK}.

In the present paper we restrict ourselves to the theory of ARS 
in Banach spaces,
which is still not quite explored. Some non-trivial examples 
of ARS in $l_2$ were found by I.Shra\u\i fel \cite{Sh l2}. 
It should be noted that
each example of a $c$-ARS in $l_2^n$ provides by Theorem \ref{ARS sm}
an example of an $\e$-net of the $n$-dimensional Euclidean ball, 
$\e=\e(c)<1$. See also \cite{Sh H} for results on ARS in
Hilbert spaces.

Some general results concerning ARS in Banach spaces
and, particularly, in uniformly smooth spaces, were obtained in \cite{V}.
There was introduced the notion of $(c,q)$-quick representing system, 
which is considerably stronger than that of ARS.

\begin{definition} \label{def QRS}
  Let $c>0$ and $q \in (0,1)$. 
  A set $D$ in a Banach space $X$ is called {\sl $(c,q)$-quick
  representing system} (or {\sl $(c,q)$-quick RS}) if for each $x \in X$ 
  there are
  scalars $(a_i)$ and elements $(x_i) \subset D$ such that
  $$
  x=\sum_{i=1}^\infty {a_i x_i}   \text{and} 
  \|a_i x_i\| \le cq^{i-1}   \text{for} i \ge 1. 
  $$
\end{definition} %QRS

\noindent It is clear that each $(c,q)$-quick RS is an ARS.
Despite of the strong restrictions in Definition \ref{def QRS}, there
exist Banach spaces $X$ in which every ARS is, in turn, a $(c_1,q)$-quick RS
for some $c_1$ and $q$. 
In \cite{V} it was proved that this happens in each super-reflexive
space $X$. 

In the present paper we generalize this result to all
B-convex Banach spaces. Suppose a space $X$ is B-convex and $Y$ is 
a subspace of $X$. We show that every $c$-ARS in $Y$ is a 
$(c_1,q)$-quick RS for some $c_1$ and $q$ depending only on $c$ and on $X$. 
This latter statement characterizes the class of B-convex Banach spaces. 
   
We characterize ARS and $(c,q)$-quick RS in uniformly smooth and 
B-convex Banach spaces via $\e$-nets of the unit balls. 
As a consequence, we have a theorem of B.~Maurey \cite{P} 
stating that the dimension of a subspace $Y$ of $\lin$ with
$Y^*$ of a good type is at most $c \log n$.

I am grateful to V.Kadets for the guidance, and
to P.Terenzi for his hospitality when I was visiting 
Politecnic Institute of Milan.

\section {Characterizations of ARS and quick RS}                  \label{Gen}
%----------------------------------------------------------------------------

Let $(x_i)_{i \in I}$ and $(y_i)_{i \in I}$ be sequences in Banach
spaces $X$ and $Y$ respectively, and let $c>0$. We call $(x_i)$ and 
$(y_i)$  {\sl $c$-equivalent} if there is a linear operator
$ T : \overline{\span}(x_i) \rightarrow \overline{\span}(y_i) $ 
which maps $x_i$ to $y_i$, and satisfies
$\|T\| \|T^{-1}\| \le c$. 

$D$ being a non-empty set, we denote the unit vectors in $l_1(D)$ 
by $e_d$, $d \in D$.

The following useful result is more or less known:  
the equivallence (i)$\Leftrightarrow$(iv) goes back to S.~Mazur 
(\cite{B}, p.~110), see also \cite{V}.

\begin{theorem}                                              \label{char ARS}
  Given a complete normalized set $D$ in a Banach space $X$, the 
  following are equivalent:

  (i) $D$ is an ARS;

  (ii) there is a $c>0$ such that each $x \in B(X)$ can be represented
  by a series $x=\sum_{i=1}^\infty {a_i x_i}$ with  
  $\sum {\|a_i x_i\|} \le c$.
  Then we call $D$ a "$c$-ARS";

  (iii) there is a quotient map $q : l_1(D) \rightarrow Z$ such that
  the sequence $(d)_{d \in D}$ is $c$-equivalent to $(qe_d)_{d \in D}$;

  (iv) there is a $c>0$ such that for every $x^* \in S(X^*)$ one has \linebreak
% WWWWWWWWWWWWWWWWWWWWWWWWWWWWWWWWWWWWWWWWWWWWWWWWWWWWWWWWWWWWWWWWWWWWWWWWWWW
% WWWWWWWWWWWWWWWWWWWWWWWWWWWWWWWWWWWWWWWWWWWWWWWWWWWWWWWWWWWWWWWWWWWWWWWWWWW
  $\sup_{d \in D} {|x^*(d)|} \ge c^{-1}$.

  In (ii), (iii) and (iv) the infimums of possible constants $c$ are
  equal and are attained.
\end{theorem}

Let us observe some nice consequences. The first one states
that ARS are stable under fairly large perturbations. 
Let $A$ and $B$ be sets in a Banach space. By definition, put
$\rho(A,B)=\sup_{a \in A} { \inf_{b \in B} {\|a-b\|} }$.

\begin{corollary}                                             \label{ARSstab}
  Let $D$ and $D_1$ be normalized sets in $X$.
  If $D$ is a $c$-ARS and $\rho(D,D_1) = \e < c^{-1}$, 
  then $D_1$ is a $c_1$-ARS, where $c_1=(1-\e c)^{-1}c$.
\end{corollary}

\proof    It follows easily from (iv) of Theorem \ref{char ARS}. 
\endproof

\begin{proposition}
  Let $D$ be a $c$-ARS in a Banach space $X$.

  (i) If $X$ is separable, then some countable subset $D_1$ of $D$
  is also a $c$-ARS.

  (ii) Let $\dim X=n$ and $c_1>c$.
  Then some subset $D_1$ of $D$ is a $c_1$-ARS 
  and $|D_1| \le e^{an}$, where 
  $a = 2 \left(  c^{-1}-c_1^{-1}  \right)^{-1}$.

  (iii) Let $\dim X=n$ and $\e>0$. Then every $x \in B(X)$ can be
  represented by a sum $x= \sum_{i=1}^n {a_i x_i}$ with 
  $(x_i) \subset D$ and $\sum {\|a_i x_i\|} \le c+\e$.
\end{proposition}

\proof     Clearly, we may assume that $D$ is normalized. Then (i)
follows in the standard way from (iv) of Theorem \ref{char ARS}.

(ii).  Let $\e=c^{-1}-c_1^{-1}$. Consider a maximal 
subset $D_1$ of $D$ such that $\|x-y\| > \e$
for $x,y \in D_1$, $x \ne y$.
By maximality, $\rho(D,D_1) \le \e$. 
Applying Corollary \ref{ARSstab}, we see that $D_1$ is a $c_1$-ARS. 
Note that the balls 
$ (d_1 + (\e/2)B(X) )_{d_1 \in D_1} $   
are mutually disjoint and are contained in 
$( 1+\e/2 )B(X)$. 
By comparing the volumes we get $|D_1| \le e^{2n/\e}$.

(iii).  By (ii), we can extract from $D$ 
a finite $(c+\e)$-ARS  $(x_i)_{i \le m}$.
By (iii) of Theorem \ref{char ARS}, 
there is a quotient map $q : l_1^m \ra Z$ such that 
the sequences $(x_i)_{i \le m}$ and $(qe_i)_{i \le m}$
are $(c+\e)$-equivalent. 
Let $T : X \ra Z$ be the isomorphism corresponding to this equivalence.
We have $\dim Z = n$ and $B(Z) = \aconv(qe_i)_{i\le m}$.
Now we use a simple consequence of Caratheodory's theorem:

\begin{itemize}
  \item   Let $K$ be a finite set in $\R^n$. Let a vector $z$ lie
  on the boundary of $\aconv(K)$. Then $z \in \aconv(z_1, \ldots, z_n)$ 
  for some $z_1, \ldots, z_n \in K$.
\end{itemize}

\noindent Applying this theorem to $K=(qe_i)_{i \le m}$, we see that 
each $z \in S(Z)$ can be represented by a sum 
$z=\sum_{k=1}^n {a_k(qe_{i_k})}$
for some subsequence $(qe_{i_k})_{k \le n}$ of $(qe_i)$
and scalars $(a_k)$ with $\sum_{k=1}^n |a_k| =1$.

Let $x \in B(X)$. Setting $z=Tx/\|Tx\|$ in the preceding observation, 
we can write
$$
Tx = \sum_{k=1}^n b_k(qe_{i_k}) = \sum_{k=1}^n b_k(Tx_{i_k}) 
\text{with}
\sum_{k=1}^n |b_k| \le \|T\|. 
$$
Thus $x = \sum_{k=1}^n b_k x_{i_k}$, and 
$$
\sum_{k=1}^n \|b_k x_{i_k}\| 
  \le  \|T^{-1}\| \sum_{k=1}^n \|b_k (qe_{i_k})\| 
  \le  \|T^{-1}\| \sum_{k=1}^n |b_k| 
  \le  \|T^{-1}\|\|T\| 
  \le  c+\e.
$$
The proof is complete.                        \endproof

\noindent{\bf Remarks.\ }  1. The estimate in (ii) is sharp by order: 
  Corollary \ref{Cor ARS tp} and Theorem \ref{haha} show that 
  any ARS in a B-convex Banach space $X$ has at least exponential
  number of terms with respect to $\dim X$.

  2. In general, one can not put $\e=0$ in (iii). Indeed, consider 
  $X=l_2^2$ and let $D$ be a countable dense subset of $S(l_2^2)$. 
  Then $D$ is a 1-ARS. However, there are points 
  $x \in S(l_2^2) \setminus \aconv(D)$ ; thus (iii) fails unless $\e>0$.

\#
Now we give a general characterization of $(c,q)$-quick RS.

\begin{theorem}                                              \label{char qRS}
  Let $D$ be a normalized set in a Banach space $X$. Suppose

  (i) $D$ is a $(c,q)$-quick RS.
  
  \noindent Then, given an $\e>0$, there are $m=m(c,q,\e)$ and 
  $b=c(1-q)^{-1}$ such that 

  (ii) the set 
  \ \ $b \cdot \bigcup \{ \aconv(D_1): D_1 \subset D, |D_1| \le m \}$ \ \ 
  is an $\e$-net of $B(X)$.

  \noindent Conversely, if $\e<1$, then (ii) implies (i) with 
  $c=b/\e$ and $q=\e^{1/m}$.
\end{theorem}

\proof       Assume (i) holds. Let $m$ be so that 
\begin{equation}
  \sum_{i>m} {cq^{i-1}} \le \e.                                 \label{cqe}
\end{equation}
Let $x \in B(X)$. For some $(x_i) \subset D$ we have 
$x=\sum_{i=1}^\infty {a_i x_i}$ 
with $|a_i| \le cq^{i-1}$. Then, by (\ref{cqe}),
$$ 
\| x- \sum_{i \le m} {a_i x_i} \| 
  = \|    \sum_{i>m}     {a_i x_i} \|   \le  \e, 
$$ 
while 
$$
\sum_{i \le m} {|a_i|} \le c(1-q)^{-1} =b.
$$ 
This proves (ii).

Conversly, assume (ii) holds. Fix an $x \in B(X)$. We shall find 
appropriate expansion $x = \sum_i{a_ix_i}$ by successive iterations. 
$S_n$ will denote the partial sum $\sum_{i \le n} {a_i x_i}$ 
(we assume $S_0=0$).

Suppose that for some $k \ge 1$ the system $(a_i)_{i \le (k-1)m}$ is 
constructed. By (ii), there are scalars 
$(a_{k,i})_{i \le m}$ and vectors
$(x_{k,i})_{i \le m} \subset D$ such that 
$|a_{k,i}| \le b$ for $i \le m$ and
\begin{equation}
  \left\|  \frac{x-S_{(k-1)m}}{ \|x-S_{(k-1)m}\| } 
         - \sum_{i \le m} {a_{k,i} x_{k,i}}  \right\|  
  \le \e.                                                         \label{e}
\end{equation}
Put $a_{(k-1)m+i}=\|x-S_{(k-1)m}\| a_{k,i}$  for $1 \le i \le m$.
Note that for each $k$
$$
  x-S_{km} = x-S_{(k-1)m} -\|x-S_{(k-1)m}\| 
  \cdot \sum_{i \le m} {a_{k,i} x_{k,i}}.
$$
Therefore, by (\ref{e}),  $\|x-S_{km}\| \le \|x-S_{(k-1)m}\| \cdot \e$.
By the inductive argument we get $\|x-S_{km}\| \le \e^k$.
Hence for $k \ge 0$ and $1 \le i \le m$,
$$
  |a_{km+i}|   =   \|x-S_{km}\| a_{k+1,i}
              \le  \e^k b
              \le  \e^{(km+i)/m-1} b
               =    \e^{-1}b \cdot (\e^{1/m})^{km+i}.
$$
Hence $ |a_i| \le \e^{-1}b (\e^{1/m})^i $   for  $i \ge 1$.  
This proves (i) with $c=b\e^{-1+1/m} \le b/\e$ and $q=\e^{1/m}$.    
\endproof

\#
Theorem \ref{char qRS} yields that, actually, the tightness
of the definition of $(c,q)$-quick RS can be substantially loosened.
Let $(b_i)$ be a scalar sequence. We say that a set $D$ in 
a Banach space $X$ is a {\em $(b_i)$-representing system}, 
if every $x \in B(X)$ admits a representation by a convergent
series $x=\sum_i a_ix_i$ with $(x_i) \subset D$ 
and $(a_i) \subset \R$, $\|a_ix_i\| \le |b_i|$ for each $i$.

\begin{corollary}
  Let $D$ be a set in a Banach space $X$ and let $\sum b_i$ be
  an absolutely convergent scalar series. Suppose
  
  (i) $D$ is a $(b_i)$-representing system.

  \noindent Then there are constants $c$ and $q$ dependent 
  only on $(b_i)$, such that

  (ii) $D$ is a (c,q)-quick representing system.

  \noindent Conversely, (ii) implies (i) with $b_i=cq^{i-1}$.
\end{corollary}

\proof        Suppose (i) holds.
Let $m$ be so that $\sum_{i>m} |b_i| \le 1/2$.
It is enough to show that (ii) of Theorem \ref{char qRS}
holds for $\e = 1/2$.
Fix $x \in B(X)$ and write its representation: 
$x = \sum_{i \ge 1} a_ix_i$ with $\|a_ix_i\| \le |b_i|$.
Then 
$$     
\| x - \sum_{i\le m} a_ix_i \| 
  =   \| \sum_{i>m} a_ix_i \|
  \le \sum_{i>m} \|a_ix_i\| 
  \le \sum_{i>m} |b_i| \le 1/2.
$$
Thus (ii) holds. The converse part is obvious.     \endproof

\#
Like ARS, quick representing systems are also stable under 
fairly large perturbations. The following analogue of Corollary \ref{ARSstab}
can easily be derived from Theorem \ref{char qRS}.

\begin{corollary}
  Let $D$ and $D_1$ be normalized sets in $X$.
  If $D$ is a $(c,q)$-quick RS and $\rho(D,D_1) = \e < (1-q)/c$,
  then $D_1$ is a $(c_1,q_1)$-quick RS,
  where $c_1$ and $q_1$ depend solely on $c$, $q$ and $\e$.
\end{corollary}

\#
Another consequence of Theorem \ref{char qRS} states 
that the cardinality of every $(c,q)$-quick RS in a 
finite-dimensional space is large. 

\begin{theorem}                                                  \label{haha}
  Let $D$ be a $(c,q)$-quick RS in a $n$-dimensional Banach space $X$.
  Then $|D| \ge e^{an}$ for some $a=a(c,q)>0$.
\end{theorem}

\noindent Before we prove this result, observe that 
there are many spaces posessing ARS of small cardinalities. 
Indeed, E.~Gluskin's construction \cite{Gl} gives us 
$n$-dimensional spaces $X_n$ and $Y_n$ having ARS of
cardinality $2n$ so that the Banach-Mazur distance 
between $X_n$ and $Y_n$ is approximately $n$.

\begin{lemma}                                                   \label{hahal}
  Let $X$ be a Banach space, $\dim X=n$,
  and $E$ be a subspace of $X$, $\dim E=m$.
  For $\e \in (0,1)$ and $b>0$, define
  $$
  U_{b,\e}(E) = b(E \cap B(X)) + \e B(X).
  $$
  Then, for some $a=a(b,\e,m)>0$,
  $$
  \vol(U_{b,\e}(E)) \le e^{-an} \vol(B(X)).
  $$
\end{lemma}

\proof    Fix a $\d>0$. Let $(z_i)_{i \le k}$ be a $\d$-net of 
$b(E \cap B(X))$; by the standard volume argument, this can
be achieved for some $k \le e^{2bm/\d}$ (see \cite{MS}, Section 2.6)
Then $(z_i)_{i \le k}$ is a $(\d+\e)$-net of $U_{b,\e}(E)$. 
Thus 
$$
\vol(U_{b,\e}(E)) 
  \le k(\d+\e)^n\vol(B(X)) 
  \le e^{2bm/\d}(\d+\e)^n\vol(B(X)).
$$
Now it is enough to pick $\d$ so that $\d+\e \le 1$.    \endproof

{\bf Proof of the Theorem \ref{haha}}.\ 
Let $\e=1/2$. Theorem \ref{char qRS} implies that 
for some $m=m(c,q)$ and $b=b(c,q)$,
$$
B(X) \subset \bigcup \{ U_{b,1/2}(E): E=\span(D_1), 
                          D_1 \subset D, |D_1| \le m \}.
$$
There are at most $\dm$ distinct members $U_{b,1/2}(E)$ in this union,
so Lemma \ref{hahal} gives us for some $a=a(b,m)$,
$$
\vol(B(X)) \le \dm e^{-an} \vol(B(X)).
$$
Hence $\dm \ge e^{an}$. The desired estimate follows easily.   \endproof

\#

Now we shall find good renormings of a space with a given ARS
or $(c,q)$-quick RS.

\begin{proposition}
  Let $D$ be a $c$-ARS in a Banach space $X$. 
  Then there is a norm $|||\cdot|||$ on $X$ which satisfies
  $\|\cdot\| \le |||\cdot||| \le c\|\cdot\|$  and such that $D$ 
  is a $1$-ARS in $(X,|||\cdot|||)$.
\end{proposition}

\proof      Set $|||x||| = \inf \sum_i {\|a_i x_i\|}$, where the infimun
is taken over all representations $x=\sum_i{a_i x_i}$ with 
$(x_i) \subset D$. Then it is enough to apply (ii) of 
Theorem \ref{char ARS}.      \endproof

\#
For $(c,q)$-quick RS, only an equivalent quasi-norm can be 
constructed.

\begin{proposition}                                            \label{renQRS}
  Let $D$ be a normalized $(c,q)$-quick RS in $X$. 
  Then there is a quasi-norm $|||\cdot|||$ on $X$ which satisfies
  $(1-q)\|\cdot\| \le |||\cdot||| \le c\|\cdot\|$  and such that 

  (i) $D \subset B(X,|||\cdot|||)$.

  (ii) $D$ is a $(1,q)$-quick RS in $(X,|||\cdot|||)$;

  (iii) the set $\cup \{tD: |t| \le c\}$ is a $q$-net
  of $B(X,|||\cdot|||)$;

\end{proposition}

\proof           For an $x \in X$, define
\begin{equation}
  |||x||| := \inf \{ \sup_{i \ge 1} {|a_i|/q^{i-1}} \},       \label{newnorm}
\end{equation}
where the infimum is taken over all sequences $(x_i) \subset D$ such that
\begin{equation}
  x = \sum_{i=1}^\infty {a_i x_i}.                               \label{repr}
\end{equation}
The homogenity of $|||\cdot|||$, (i) and (ii) follow easily.

Now we show that $1-q \le |||x||| \le c$ for every $x \in S(X)$. 
The right hand side follows from (\ref{newnorm}). 
Conversely, let (\ref{repr}) be a representation of $x$ such that 
$\sup_i {|a_i|/q^{i-1}} = \lambda < \infty$. Then
$$
  1  =    \big\| \sum_{i=1}^\infty {a_i x_i} \big\|
     \le  \sum_{i=1}^\infty {|a_i|}
     \le  \sum_{i=1}^\infty {\lambda q^{i-1}}
     =    \lambda (1-q)^{-1}.
$$
Thus $\lambda \ge 1-q$; therefore $|||x||| \ge 1-q$.

It remains to prove (iii). Pick any $x \in X$ with $|||x||| \le 1$ 
and $\e>0$.
Let (\ref{repr}) be any expansion with  $|a_i|/q^{i-1} \le 1+\e$ 
for $i \ge 1$. Write
$$
  x-a_1 x_1  =  \sum_{i=1}^\infty {a_{i+1} x_{i+1}}.
$$
Then $|||x-a_1 x_1||| \le \sup_i {|a_{i+1}|/q^{i-1}} \le (1+\e)q$. 
This proves (iii).     \endproof

\noindent {\bf Remarks.\ } 
1. The statement (iii) of Proposition \ref{renQRS} means that 
in the new norm one can take $\e=q$, $b=c$ and $m=1$ 
in Theorem \ref{char qRS} (ii).

2. In general, there is no equivalent 
{\sl norm} $|||\cdot|||$
satisfying (ii) or (iii) of Proposition \ref{renQRS}. Indeed, take 
$X=l_2^2$ and $D=\{ (1,0), (0,1) \}$. Then $D$ is a $(4,1/4)$-quick RS,
but $D$ cannot be $(1,1/4)$-quick RS in any norm $|||\cdot|||$ on $X$,
nor can the set $\cup \{tD: t \in \R\}$ be a $1/4$-net of $B(X,|||\cdot|||)$.

\section {Absolutely representing systems
          in uniformly smooth spaces}                          \label{Smooth}
%----------------------------------------------------------------------------

We recall the notion of uniform smoothness (see \ \cite{DGZ}).
Let $X$ be a Banach space. The {\sl modulus of smoothness} of $X$ is 
the function defined for $\tau>0$ by
$$
\rho(\tau) = \sup\{ (\|x+y\|+\|x-y\|)/2-1: 
                         x,y \in X, \|x\|=1, \|y\| \le \tau \}.
$$
$X$ is called {\sl uniformly smooth} if \
$\lim_{\tau\rightarrow 0} {\rho(\tau)/\tau} =0$.

\begin{theorem}                                         \label{ARS sm}
  Let $D$ be a normalized set in a Banach space $X$ and $c>1$.
  Suppose $\rho(\tau)/\tau \le (4c)^{-1}$ for some $\tau \in (0,1)$. 
  Suppose 

  (i) $D$ is a $c$-ARS in $X$.
  
  \noindent Then letting $t=2\tau/3$ and $\e=1-\tau/3c$, we have
  
  (ii) the set $\pm tD$ is an $\e$-net of $B(X)$.
  
  \noindent Conversely, if $\e<1$, then (ii) implies (i) with 
  $c=c(t,\e)$.
\end{theorem}

\noindent{\bf Remark.\ } The converse part of Theorem \ref{ARS sm} 
holds in every Banach space $X$. Indeed, it is enough to apply 
Theorem \ref{char qRS} and note that each $(c,q)$-quick RS is 
a $c_1$-ARS for $c_1=c(1-q)^{-1}$.

An immediate consequence follows:

\begin{corollary}                                          \label{cor ARS sm}
  Let $D$ be a normalized $c$-ARS in a uniformly smooth 
  space $X$. Then there are constants $t>0$ and $\e<1$ 
  depending solely on $c$ and on the modulus of smoothness of $X$
  so that the set $\pm tD$ is an $\e$-net of $B(X)$.
\end{corollary}

Recall that each superreflexive space $X$ has an equivalent norm
$|||\cdot|||$ such that $(X,|||\cdot|||)$ is a uniformly smooth
space (see \cite{DGZ}).
Therefore, for each super-reflexive space $X$ the conclusion of 
Corollary \ref{cor ARS sm} will be true after an equivalent renorming. 

Moreover, this property characterizes the class of super-reflexive
spaces. Indeed, let $X$ be not super-reflexive; then $X$ is not
super-reflexive in any equivalent norm. Let $\delta>0$.
Then there are almost square sections of $B(X)$ (see \cite{DGZ}).
More precisely, there is a system of two vectors
$(z_1,z_2)$ in $S(X)$ which is $(1+\delta)$-equivalent to the canonical
vector basis of $l_\infty^2$. Let $Z=\span(z_1,z_2)$. Then $Z$ is 
$(1+\delta)$-isomorphic to $l_\infty^2$ and hence is 
$(1+\delta)$-complemented in $X$; write $X=Z \oplus Y$ for an 
corresponding complement $Y$ in $X$. Put $D=\{z_i+y: y \in Y, \; i=1,2\}$.
Now it is not hard to check that $D$ is a $3$-ARS in $X$, but 
the set $\cup\{tD: t \in \R\}$ is not an $\e$-net of $B(X)$ unless
$\e > 1-\delta/2$. 
This argument was shown to me by V.~Kadets.

\#
The proof of Theorem \ref{ARS sm}
requires some $(\e<1)$-net tools.

\begin{lemma}                                                      \label{L9}
  Let $\lambda \in [0,1]$ and $A\subset \lambda \cdot B(X)$.
  Suppose that $A$ is a $\lambda $-net for $S(X)$. 
  Then         $A$ is a $\lambda $-net for $B(X)$.
\end{lemma}

\proof       For each $x\in B(X)$, there
exists an $y\in A$ such that $\|x/\|x\|-y\|\leq \lambda $.
Hence 
\begin{eqnarray*}
  \|x-y\|   &=  &   \big\|  \|x\|(x/\|x\|-y)-(1-\|x\|)y  \big\|   \\
            &\le&   \|x\|\lambda +(1-\|x\|)\lambda =\lambda . 
\end{eqnarray*}
This completes the proof.          \endproof

\begin{lemma}                                                     \label{L10}
  Let $A\subset X$ be a $(1-\delta)$-net for $S(X)$ with
  $\delta \in (0,1)$. Then, for each $\gamma \in [0,1]$, the set 
  $\gamma A$ is a $(1-\gamma \delta )$-net for $S(X)$.
\end{lemma}

\proof      For any $x\in S(X)$ there exists an $y\in A$ such
that $\|x-y\| \le 1-\delta$. Hence 
\begin{eqnarray*}
  \|x-\gamma y\|   & = &   \| \gamma(x-y)+(1-\gamma)x \|    \\
                   &\le&   \gamma(1-\delta)+(1-\gamma)      
                   = 1-\gamma \delta,
\end{eqnarray*}
which concludes the proof.            \endproof

\begin{corollary}                                                 \label{C11}
  Let $\tau >0$, $\delta \in (0,1)$ and let 
  $A\subset \tau \cdot B(X)$ be a $(1-\delta)$-net for $S(X)$. 
  Then, for each $ 0 \le \gamma \le \min (1,\frac 1{\tau +\delta})$, 
  the set $\gamma A$ is a $(1-\gamma \delta)$-net for $B(X)$.
\end{corollary}

\proof          By Lemma \ref{L10}, $\gamma A$ is a 
$(1-\gamma \delta)$-net for $S(X)$. On the other hand, 
$\gamma \tau \le 1-\gamma \delta$, so that 
$\gamma A \subset (1-\gamma \delta) \cdot B(X)$.
Then, by Lemma \ref{L9}, $\gamma A$ is a $(1-\gamma \delta)$-net 
for $B(X)$.            \endproof

\#
Now, we establish a "locally equivalent norm" on $X$.

\begin{lemma}                                                     \label{L12}
  Let $x\in S(X)$ and $x^* \in S(X^*)$ be such that $x^*(x)=1$.  
  Then for each $z\in X$ we have: 
  $$
  x^*(z) \le \|z\| \le x^*(z)+2\rho(\|z-x\|). 
  $$
\end{lemma}

\proof           Put $y=x-z$. Then 
\begin{eqnarray*}
  2\rho(\|y\|)   &\ge&   \|x+y\|+\|x-y\|-2               \\
                 &\ge&   x^*(x+y)+\|x-y\|-2              \\
                 &\ge&   1+x^*(y)+\|x-y\|-2              \\
                 & = &   \|x-y\|-x^*(x-y) 
                   = \|z\|-x^*(z).
\end{eqnarray*}
Hence the right inequality is proved while the left one is trivial.
\endproof

\#
{\bf Proof of the Theorem \ref{ARS sm}}.\ 
Assume (i) holds. We claim that the set $\pm \tau D$ is a 
$(1-\tau/2c)$-net of $S(X)$. 
Indeed, given an $x\in S(X)$, one can pick a functional 
$x^*\in S(X^*)$ such that $x^*(x)=1$. Then, 
by Theorem \ref{char ARS}, we have 
$$
\theta x^*(x) \ge c^{-1} 
$$
for some $x \in D$ and some $\theta \in \{-1,1\}$. 
Now apply Lemma \ref{L12} with $z=x-\theta \tau x$ :
\begin{eqnarray*}
  \| x-\theta \tau x \|   &\le&   x^*(x-\theta \tau x)+2\rho(\tau )   \\
                          &\le&   1-\tau c^{-1} + 2\rho(\tau )        \\
                          &\le&   1-\tau c^{-1} + 2\cdot \tau/4c
                            = 1-\tau/2c.
\end{eqnarray*}
This proves our claim.

Then apply Corollary \ref{C11}: $A=\pm \tau D$, $\delta=\tau/2c$ 
and $\gamma=2/3$ will satisfy its conditions. We get that
$\frac{2}{3} A$ turns to be a $(1-\tau/3c)$-net of $B(X)$,
proving (ii). 

The converse part follows from the remark above.     \endproof

\section {Absolutely representing systems 
          and type of Banach spaces}                             \label{Type}
%----------------------------------------------------------------------------

The theory of type and cotype for normed spaces 
can be found in \cite{MS} or \cite{LeTa}.
By $(\e_i)$ we denote a sequence of independent random 
variables with the distibution $\P\{\e_i=1\} = \P\{\e_i=-1\}=1/2$.
Consider a Banach space $X$ of type $p>1$, i.e. such that there is a $c>0$
such that the inequality
\begin{equation}
  \E \big\| \sum_{i \le n} {\e_i x_i} \big\| ^p 
     \le c^p \sum_{i \le n} {\|x_i\|^p}                             \label{t}
\end{equation}
holds for each $n>0$ and each sequence $(x_i)_{i \le n}$ in $X$.
By $\tpx$ we denote the least constant $c$ for which 
the inequality (\ref{t}) always holds. For $p>1$, we denote by $p^*$ 
the conjugate number: $1/p+1/p^*=1$. 

The following result is contained implicitly in \cite{P}
and is known as a "dimension-free variant of Caratheodory's theorem".
For the sake of completeness, we include its proof.

\begin{theorem}                                                \label{ARS tp}
  Let $D$ be a normalized set in a Banach space $X$ of type $p>1$.
  Suppose that for some $c>1$

  (i) $D$ is a $c$-ARS.

  \noindent Let $k>0$. Put $c_1=c$ and $\e = 4c \tpx k^{-1/p^*}$.
  Then
  
  (ii) the set $\{ c_1 k^{-1} \sum_{i \le k} {\pm x_i} : (x_i) \subset D \}$
  is an $\e$-net of $B(X)$.
 
  \noindent Conversely, (ii) implies (i) with $c=c(c_1,k)$.
\end{theorem}

\noindent Applying Theorem \ref{char qRS}, we obtain

\begin{corollary}                                          \label{Cor ARS tp}
  Let $D$ be a normalized set in a Banach space $X$ of type $p>1$.
  Suppose that for some $c>1$

  (i) $D$ is a $c$-ARS.

  \noindent Then, for some $c_1=c_1(c,p,\tpx)$ and $q=q(c,p,\tpx)$, we have:
  
  (ii) $D$ is a $(c_1, q)$-quick RS.
  
  \noindent Conversely, (ii) implies (i) with $c=c(c_1,q)$.
\end{corollary}

\noindent Before the proof of Theorem \ref{ARS tp}, let us give 
some comments.
A Banach space $X$ is called {\sl B-convex} if it does not contain
$l_1^n$ uniformly. $X$ is B-convex iff $X$ is of some type $p>1$. 
It follows that if $X$ is a B-convex Banach space and 
$D$ is a $c$-ARS in some subspace of $X$, 
then $D$ is a $(c_1,q)$-quick RS, where the constants $c_1$ and $q$
depend only on $c$ and $X$.

Moreover, the latter property characterizes B-convex Banach spaces. 
Indeed, fix a space $X$ which is not B-convex. Then, for 
each positive integer $n$, there is a 
sequence $(x_{n,i})_{i \le n}$ in $X$ which is $2$-equivalent to
the canonical vector basis of $l_1^n$. Take $D_n=(x_{n,i})_{i \le n}$
and $Y_n = \span(D_n)$. Then $D_n$ is a 2-ARS in $Y_n$. 
However, letting $n \rightarrow \infty$,
we see that $D_n$ cannot be a $(c_1, q)$-quick RS for fixed $c_1$ and $q$.

One exciting problem remains unsolved. We have got that each
ARS in a B-convex space $X$ is a $(c,q)$-quick RS for some $c$ and $q$.
Does this happen only in B-convex spaces?

\#
{\bf Proof of Theorem \ref{ARS tp}}.\ 
Fix any $x \in B(X)$. Then, for some $(x_i) \subset D$, there is
a representation $x=\sum_{i=1}^\infty {a_ix_i}$ with $\sum|a_i| \le c$.

Then there is a sequence $(\xj)_{j \ge 1}$ of independent random variables
with the following distribution for every $i,j \ge 1$:
\begin{eqnarray*}
  & &\P \{ \xj = \sign (a_i) c x_i \}   =   c^{-1} |a_i|,        \\
  & &\P \{ \xj = 0 \}   =   1-c^{-1} \sum_{n}|a_n|.
\end{eqnarray*}
Therefore $\E\xj=x$ for each $j$. Now, since $\xj$ are independent, 
we have
$$
\E \big\| \sum_{j \le k} {(\xj - \E\xj)} \big\|^p
    \le (2\tpx)^p \sum_{j\le k} {\E\|\xj-\E\xj\|^p}
$$
(see \cite{LeTa}, Chapter 9). Note that $\E\|\xj-\E\xj\|^p \le (c+1)^p$; 
hence
$$
\E \big\| k^{-1}\sum_{j \le k} {(\xj - \E\xj)} \big\|^p
    \le (2\tpx)^p k^{-p} \cdot k(c+1)^p.
$$
Therefore
$$
\E \big\| -x+k^{-1}\sum_{j \le k} \xj \big\|^p
    \le \left( 2\tpx (c+1)k^{-1/p^*} \right)^p.
$$
In particular, there is one realization of the random variable
$(-x+k^{-1}\sum_{j \le k} \xj)$   so that  
$$
\big\| -x+k^{-1}\sum_{j \le k} \xj \big\|
    \le 2\tpx (c+1)k^{-1/p^*}.
$$                                         
This concludes the proof.      \endproof

\#
In conclusion, let us show how these results provide an estimate
from above on the dimension of nice sections of the cube.
The following result due to B.~Maurey is proved in \cite{P}.

\begin{theorem} (B.~Maurey). 
  Let $X$ be a finite dimensional space, $p>1$ and $\tpxx \le C$. 
  Suppose that $X$ is $c$-isomorphic to some subspace of $l_\infty^n$. 
  Then, for some $a=a(p,C,c)$, we have
  $$
  \dim X \le a\log n.
  $$
\end{theorem}

\proof   By duality, $X^*$ is $c$-isomorphic to some quotient space
of $\l_1^n$. Then, Theorem \ref{char ARS} gives us a $c$-ARS 
$D$ in $X^*$ with $|D|=n$. By Corollary \ref{Cor ARS tp}, 
$D$ is a $(c_1,q)$-quick RS in $X^*$ 
for some $c_1=c_1(p,C,c)$ and $q=q(p,C,c)$.
Then Theorem \ref{haha} yields $n \ge e^{a \dim X}$
for some $a=a(c_1,q)>0$.    \endproof

\#

{\em E-mail: } paoter@mate.polimi.it

\end {document}